%% file: Weilnumbers.tex
\documentclass[12pt]{article}
\usepackage{amssymb,amsmath,amscd}
\usepackage[all]{xy}
\input{def_engl.tex}

\begin{document}
\title{On Weil numbers in cyclotomic fields\footnote{MSC Classification: 11R18, 11R23,
11L05}}
\author{Bruno Angl\`es, Tatiana Beliaeva}

\maketitle

Let $p$ be an odd prime number. It was noticed by Iwasawa that the
$p$-adic behavior of Jacobi sums in $\qq(\z_p)$ is linked to
Vandiver's Conjecture (see \cite{Iwasawa}). This result has been
generalized by various authors for the cyclotomic $\zp$-extensions
of abelian fields (see for example \cite{Hachimori_Ichimura},
\cite{Ichimura}, \cite{moi}). In this paper we consider the module
of Weil numbers (see \S 2) for the cyclotomic $\zp$-extension of
$\qq(\z_p)$, and we get some results quite similar to those for
Jacobi sums. In particular we establish a connection between the
$p$-adic behavior of Weil numbers and a weak form of Greenberg
Conjecture (see \cite{Thong2}, \cite{Raphi_Thong})

\section {Notations}
Let $p$ be a fixed odd prime number. For any $n\in \nn$ we denote
by $k_n$ the $p^{n+1}$-th cyclotomic field $\qq(\mu_{p^{n+1}})$,
where $\mu_{p^{n+1}}$ is the group of $p^{n+1}$-th roots of unity.
We note $\D=\Gal{(k_0/\qq)}$, $\Ga_n=\Gal({k_n/k_0})$ and
$G_n=\Gal{(k_n/\qq)}$, so $G_n\simeq\D\times \Ga_n$. Let $\z_p\in
\mu_p\setminus\{1\}$ and take for any $n\in\nn$ $\z_{p^{n+1}}\in
\mu_{p^{n+1}}$ such that $\forall n\geq 1$
$\z_{p^{n+1}}^p=\z_{p^n}$. We note $\pi_n=1-\z_p^{n+1}$

We shall also use the following more or less standard notations:

\noindent $k_{n,p}$ the $p$-completion of $k_n$;

\noindent $\U_n=1+\pi_n\zp[\z_{p^{n+1}}]$ principal units in
$k_{n,p}$;

\noindent $\Ga=\limp\Ga_n\simeq\zp$, $\g_0$ its topological
generator, where $\forall \eps\in\mu_{p^\infty}$,
$\g_0(\eps)=\eps^{1+p}$;

\noindent $\La=\zp[[\Ga]]$ the Iwasawa algebra of the profinite
group $\Ga$, $\La\simeq \zp[[T]]$ by sending $\g_0-1$ to $T$
(\cite[Theorem 7.1]{Washington});

\noindent $A_n$ is the Sylow $p$-subgroup of $Cl(k_n)$, where
$Cl(k_n)$ is the ideal class group of $k_n$;

\noindent $X=\limp A_n$ be the projective limit of $A_n$ for the
norm maps;

\noindent $I_n$ the group of prime-to-$p$ ideals of $k_n$;

\noindent $k_\infty=\bigcup_{n\in \nn}k_n$,
$\Gal(k_\infty/k_0)=\Ga$;

\noindent $L_n/k_n$ the maximal abelian unramified $p$-extension
of $k_n$; $\Gal(L_n/k_n)\simeq A_n$ by class field theory;

\noindent $M_n/k_n$ the maximal abelian $p$-extension of $k_n$
 unramified outside of $p$;

\noindent $\X_n=\Gal(M_n/k_n)$;

\noindent $L_\infty=\bigcup L_n$; $X\simeq
\Gal(L_\infty/K_\infty)$;

\noindent $M_\infty=\bigcup M_n$;

\noindent $\X_\infty=\Gal(M_\infty/k_\infty)$.

Let $\psi$ be a fixed odd character of $\D$, different from
Teichm\"uller character $\om$. We note $e_\psi$ the associated
idempotent defined by
$$
e_\psi=\frac{1}{|\D|}\sum_{\delta\in\D}\psi(\delta)\delta^{-1}\in
\zp[\D].
$$
Let $\Mc\in\La$ be the distinguished polynomial of smallest degree
such that $\Mc(T) e_\psi X=\{0\}$. We call it  \emph{the minimal
polynomial of $e_\psi X$.} It is well known to be prime to
$\omega_n=(T+1)^{p^n}-1$ for any $n$ (cf. \cite[\S 13.6, Theorem
7.10, Theorem 5.11 and Theorem 4.17]{Washington}).

\section {Weil numbers and Jacobi sums.}

Fix an $n$ for a moment.

\begin{Def} We call \emph{Weil module} of $k_n$ the module $\W_n$ defined by

\begin{multline}\label{def_weil}
\W_n=\left\{ f\in \Hom_{\zz[G_n]}(I_n, k_n^*)\,\mid\
\exists\beta(f)\in \zz[G_n] \text{ such that } \forall \ag\in I_n
\right. \\ \left.\ag=(\alpha)\Rightarrow f(\ag)\equiv
\alpha^{\beta(f)} \mod \mu_{2p^{n+1}} \right\}
\end{multline}
\end{Def}
Observe that $\forall f\in \W_n$, $f(I_n)\subset
\mu_{2p^{n+1}}\U_n$.
\begin{Def}
So we define the {\rm module of Weil numbers} ${W}_n$
$$
{W}_n=\{ f(\ag)\,|\, f\in \W_n,\ \ag\in I_n\}.
$$
Observe that $W_n$ is a submodule of $\mu_{2p^{n+1}}\U_n$.
\end{Def}

Let $k_n^+$ be the maximal totally real subfield of $k_n$ and let
$G_n^+$ stay for $\Gal(k_n^+/\qq)$. Let $N_n$ be the norm element
in $\zz[G_n]$. Let $N_n^+\in \zz[G_n]$ be such that its image by
the restriction map $\zz[G_n]\longrightarrow\zz[G_n^+]$ is
$\sum_{\sigma\in G_n^+}\sigma$.

\begin{Lem}
Let $f\in \W_n$. Then $\beta(f)\in\zz[G_n]$ is unique and
$$ \beta(f)\in N_n^+\zz + \zz[G_n]^-.$$
\end{Lem}

\noindent {\bf Proof:} Let $f$ be in $\W_n$ and suppose $\beta(f)$
and $\beta'(f)$ verify the required condition.

Let $\p$ be a split prime ideal in $I_n$. Let $m\geq 1$ be such
that $\p^m=\alpha\Oc_{k_n}$. Then
$$
f(\p^m)\equiv\alpha^{\beta(f)}\equiv\alpha^{\beta'(f)}\mod\mu_{2p^{n+1}}.
$$
 Thus
$\p^{m\beta(f)}=\p^{m\beta(f')}$, that implies
$\p^{m(\beta(f)-\beta(f'))}=\Oc_{k_n}$, so $\beta(f)=\beta'(f)$.
Furthermore:
$$
\beta(f)\in
\Ann_{\zz[G_n]}(\Oc_{k_n}^*/\mu_{2p^{n+1}})=N_n^+\zz+\zz[G_n]^-. \
\square
$$
\medskip

\begin{Prop}\label{prop1}
The map $\beta\,:\,\W_n \longrightarrow N_n^+\zz+\zz[G_n]^-$
defined by $f\mapsto \beta(f)$ gives rise to the exact sequence of
$\zz[G_n]$-modules.
$$
0\longrightarrow\Hom_{\zz[G_n]}(I_n,\mu_{2p^{n+1}})\longrightarrow
\W_n^- \longrightarrow (\Ann_{\zz[G_n]}Cl(k_n))^-\longrightarrow
B_n\longrightarrow0
$$
where $B_n$ is a finite abelian elementary $2$-group.
\end{Prop}

\noindent {\bf Proof:}

By the definition of $\W_n$ one has
$$
\Hom_{\zz[G_n]}(I_n,\mu_{2p^{n+1}})=\{f\in \W_n\,|\,
\beta(f)=0\}=\ker \beta.
$$

Note that $f\in\W_n^-$ implies $\beta(f)\in \zz[G_n]^-$. Take
$f\in \W_n^-$ and $\p$ a prime ideal in $I_n$. Let $\p$ be a split
prime ideal in $I_n$. Let $m\geq 1$ be such that $\p^m$ is
principal. Then
$$
\p^{m\beta(f)}= f(\p)^{m}\Oc_{k_n},
$$
that implies
$$
\p^{\beta(f)}=f(\p)\Oc_{k_n}.
$$
Thus $\beta(f)\in(\Ann_{\zz[G_n]}Cl(k_n))^-$.

Let $\beta$ be in $(\Ann_{[\zz[G_n]}Cl(k_n))^-$ and $\p$ a prime
ideal in $I_n$. Then there exists $\g\in k_n^*$ such that
$\p^\beta=\g\Oc_{k_n}$. Let $\bar{\g}$ be the complex conjugate of
$\g$. Then $\bar{\g}=\g^{-1}\eps$ for some $\eps\in \Oc_{k_n}^*$.
Thus $\eps=\g\bar{\g}$, i.e. $\eps$ is a real unit. Consider
$\g_1=\eps^{-1}\g^2$. One has: $\g_1^{-1}=\bar{\g_1}$ and
$$
\p^{2\beta}=\g^2\Oc_{k_n}=\eps^{-1}\g^2\Oc_{k_n}=\g_1\Oc_{k_n}.
$$
Let $\g_2\in k_n^*$ such that
$\p^{2\beta}=\g_1\Oc_{k_n}=\g_2\Oc_{k_n}$ and
$\bar\g_2=\g_2^{-1}$. Then $\g_1=\g_2\eta$ for some $\eta\in
\Oc_{k_n}^*$, $\g_1^{-1}=\g_2^{-1}\bar\eta$. That implies
$\eta\bar\eta=1$, i.e. $\eta$ is a root of unity.

Now one can choose, for any $\p\in I_n$, $\g_\p\in k_n^*$ such
that $\p^{2\beta}=\g_\p\Oc_{k_n}$, $\bar\g_\p=\g_\p^{-1}$ and
$\g_{\p^\sigma}=\g_\p^\sigma$ $\forall\sigma\in G_n$. We set:
$$
f(\p)=\g_\p
$$
and one can verify that $f\in \W_n^-$ and $\beta(f)=2\beta$. Thus
$$
2(\Ann_{[\zz[G_n]}Cl(k_n))^-\subset\beta(\W_n^-)\subset(\Ann_{[\zz[G_n]}Cl(k_n))^-,
$$
that completes the proof. $\square$

\bigskip

Let $l\neq p$ be a prime number. Let $\lc$ be the prime ideal of
$k_n$ above $l$ and $q=\mid \Oc_{k_n}/\lc\mid$. Fix a primitive
$l$-th root of unity $\z_l$. The Gauss sum $\tau_n(\lc)$
associated to $\lc$ is defined by
$$
\tau_n(\lc)=-\sum_{a\in{\mathbb
F}_q}\chi_\lc(a)\z_l^{\Tr_{{\mathbb
F}_q/{\mathbb F}_l}(a)}
$$
where $\chi_\lc$ is a character on $\mathbb F_q^* $ of order
$p^{n+1}$ defined by
$$
\chi_\lc(x)\equiv x^{-\frac{q-1}{p^{n+1}}}\mod \lc.
$$
One can show that $\forall\delta\in G_n$ on has
$\tau_n(\lc^\delta)=\tau_n(\lc)^\delta$ (see \cite[\S
6.1]{Washington}. So we have a well defined morphism of
$\zz[G_n]$-modules
\begin{equation*}
    \tau_n\,:\, I_n\longrightarrow \Omega(\z_{p^{n+1}})^*,
\end{equation*}
where $\Omega$ is the compositum of all the $\qq(\z_m)$, $m$ prime
to $p$.
\begin{Def} The \emph{Jacobi module} $\Jc_n$  associated to $k_n$
is defined by
$$
\Jc_n=\zz[G_n]\tau_n\cap\Hom_{\zz[G_n]}(I_n,k_n^*),
$$
and the \emph{module of Jacobi sums} $J_n$ is defined by
$$
J_n=\{f(\ag)\, \mid\, f\in \Jc_n, \ \ag\in I_n\}.
$$
\end{Def}

Let us denote by $\sigma_a$ the image of $a\in \zz$, prime to $p$,
via the standard isomorphism $(\zz/p^{n+1}\zz)^*\simeq G_n$. Let
$$
\theta_n=\frac{1}{p^{n+1}}\sum_{a=1,\
(a,p)=1}^{p^{n+1}}a\sigma_a^{-1}
$$
be the Stickelberger element of $k_n$. Set
$$
{\cal S}_n'=\sum_{(t,p)=1}\zz[G_n](t-\sigma_t).
$$
\begin{Def}
The \emph{Stickelberger ideal of $k_n$} is defined by
$$
{\cal S}_n={\cal S}_n'\theta_n,
$$
(see \cite[Lemma 6.9]{Washington}).
\end{Def}

\begin{Th}[Stickelberger's Theorem {\cite[Theorem3.1]{Sinnott}}] Let $\p$ be a prime ideal in $I_n$, and
$\beta\in{\cal S}_n'$. Then $\tau(\p)^\beta\in k_n^*$,
$\beta\theta_n\in\zz[G_n]$ and
$$
\tau_n(\p)^\beta\Oc_{k_n}=\p^{\beta\theta_n}.
$$
Moreover,  $\tau_n^\beta(\p) \in \U_n$.
\end{Th}
\begin{Lem}\label{S'}
$$
{\cal S}_n'=\{\beta\in\zz[G_n]\,|\, \tau_n^\beta\in \Jc_n\}.
$$
\end{Lem}
\noindent{\bf Proof:} The inclusion ${\cal
S}_n'\subset\{\beta\in\zz[G_n]\,|\, \tau^\beta\in \Jc_n\}$ is
obvious by Stickelberger's theorem.

To prove the inverse inclusion it suffices to show that
$\tau^\beta\in \Jc_n$ implies $\beta\theta_n\in \zz[G_n]$.

Let $\p\in I_n$ be a split prime ideal and $\widetilde\p$ the
unique prime ideal of $\zz[\z_{p^{n+1}},\z_l]$ above $\p$, where
$l=\p\bigcap\qq$, $l\equiv 1\mod p^{n+1}$. Then
$$
\tau_n(\p)\zz[\z_{p^{n+1}},\z_l]=\widetilde\p ^{(l-1)\theta_n}.
$$
Thus
$$
\tau_n(\p)^\beta\zz[\z_{p^{n+1}},\z_l]=\widetilde\p
^{(l-1)\beta\theta_n}.
$$
On the other hand,
$$
\tau_n(\p)^\beta\Oc_{k_n}=\p^z
$$
for some $z\in \zz[G_n]$, that implies
$$
\widetilde{\p}^{(l-1)z}=\widetilde{\p}^{(l-1)\theta_n\beta}
$$
and since $l\equiv1\mod(p^{n+1})$, one has
$(l-1)z=(l-1)\theta_n\beta$. Thus $z=\theta_n\beta$ that implies
$\beta\theta_n\in\zz[G_n]$. $\square$

\begin{Prop} \label{prop2}~
\begin{itemize}
\item[(1)] $\Jc_n\subset \W_n$

\item[(2)] $\Jc_n\simeq {\cal S}_n$.
\end{itemize}
\end{Prop}
\noindent{\bf Proof:}

\noindent(1)Using the lemma \ref{S'} one can easily verify that
$$
\Jc_n=\{\tau_n^\delta\,| \, \delta \in {\cal S}_n'\}.
$$
Then for any $f\in \Jc_n$ there exists $\delta \in {\cal S}_n'$
such that $f=\tau_n^\delta$.

 Let $f\in \Jc_n$ and let $\ag\in I_n$ be a  principal ideal,
$\ag=\alpha\Oc_{k_n}$. Then by the Stickelberger Theorem one has
$$
f(\ag)=\tau_n^\delta(\ag)=\eps\alpha^{\delta\theta_n}
$$
for some unit $\eps$. But
$$
\tau_n(\ag)\overline{\tau_n(\ag)}=N_n(\ag)=N_n(\alpha),
$$
so $\eps\in\mu_{2p^{n+1}}$. That means $f(\ag)\equiv
\alpha^{\delta\theta_n} \mod \mu_{2p^{n+1}}$, i.e. $f\in\W_n$.

\medskip
\noindent (2) As $\Jc_n\subset \W_n$ by (1), the map
$\beta|_{\Jc_n}$ is well defined. On the other hand, for any $f\in
\Jc_n$
$$
\beta(f)=\beta(\tau_n^\delta)=\delta\theta_n
$$
for some $\delta \in {\cal S}_n'$. Thus one has a well defined map
$$
\begin{array}{ccc}
   \Jc_n & \longrightarrow & {\cal S}_n \\
   \tau_n^\delta & \longmapsto & \delta\theta_n \\
\end{array}
$$
This map is obviously surjective (by the definition of ${\cal
S}_n$). Its kernel is a submodule of
$\Hom_{\zz[G_n]}(I_n,\mu_{p^{n+1}})$ by proposition \ref{prop1}.
Let $\delta\in \zz[G_n]$ such that $\delta\theta_n=0$. Then we
have $\sigma_{-1}\delta=\delta$ ($\sigma_{-1}$ being the complex
conjugation in $G_n$) and $\delta N_n=0$.

Now let $\delta\in{\cal S'}_n $ such that $\delta\theta_n=0$. Then
$$\tau_n^{\sigma_{-1}\delta}=\tau_n^\delta,$$
 and
$$
\tau_n^{\delta\sigma_{-1}}=(\tau_n^{\sigma_{-1}})^\delta=\tau_n^{-\delta}
$$
as $\theta_n\theta_n^{\sigma_{-1}}=N_n$. Thus
$\tau_n^{2\delta}=1$. Therefore $\tau_n^\delta=1$ as
$\tau_n^\delta\equiv1\mod \pi_n$.
 $\square$


\medskip
\begin{Lem} \label{gauss-norms} Let
$N_{n,n-1}$ be the norm map in the extension $k_n/k_{n-1}$ and
$\Lc\in I_n $ a prime ideal. Then
$$
N_{n,n-1}(\tau_n(\Lc))=\tau_{n-1}( N_{n,n-1}(\Lc))\zeta^a l^b,
$$
for some $a,b\in\zz$ and some $\z\in \mu_{p^{n+1}}$.
\end{Lem}
For a proof see \cite[Lemma 2]{Ichimura}.
\begin{Rem}
The composition $N_{n,n-1}\circ \tau_n$ is well defined because\\
$\Gal(k_n/k_{n-1})\simeq\Gal(\Omega(\z_{p^{n+1}}))/\Omega(\z_{p^n}))$,
$\forall \geq 1$.
\end{Rem}
\begin{Lem}[{\cite[Proposition 7.6 (c)]{Washington}}] \label{lem_S}
The restriction map $\Res:\zz[G_n]\rightarrow \zz[G_{n-1}]$
induces the surjective map
$$
\Res\,:\, {\cal S}_n^-\longrightarrow {\cal S}_{n-1}^-.
$$
\end{Lem}

\begin{Prop}~
$$
\forall n\geq1,\ N_{n,n-1}\circ \Jc_n^-\equiv J_{n-1}^-\circ
N_{n,n-1}\mod \Hom_{\zz[G_n]}(I_n,\mu_{2p^{n+1}})
$$
\end{Prop}
\noindent{\bf Proof:} Let $f\in {\cal J}^-_n$. By the proposition
\ref{prop2} there exists some $\beta\in {\cal S}_n^-$ such that
$f=\tau_n^\beta$. Let $\Lc\in I_n$ be a prime ideal and
$\lc=N_{n,n-1}$. Then by the lemmas \ref{gauss-norms} and
\ref{lem_S}
$$
N_{n,n-1}(\tau_n^\beta(\Lc))\equiv\tau_{n-1}^{\Res(\beta)}(\lc)
\mod \mu_{p^{n+1}}(\Lc).
$$
The Proposition follows. $\square$

\section{Annihilators}
We recall that $\psi$ is an odd $\qp$-valued character of $\D$,
irreducible over $\qp$, different from Teichm\"uller character
$\om$.
\begin{Lem} Let $\Mc\in \La$ be the minimal polynomial of $e_\psi X$. Then
$$
\limp e_\psi(\Ann_{\zp[G_n]}A_n)=\Mc(T)\La,
$$
the projective limit being taken for the restriction maps.
\end{Lem}
\noindent{\bf Proof:} First we remark that
$$
e_\psi\Ann_{\zp[G_n]}A_n=\Ann_{\zp[\Ga_n]}e_\psi A_n.
$$
We set $A_{n,\psi}=e_\psi A_n$ for simplicity.

Let $\Mc=(\Mc_n)_{n\geq0}\in\La\simeq \limp\zp[\Ga_n]$, the limit
being taken with respect for restriction maps. As $X^-$ has no
nontrivial finite submodule (see \cite[Proposition
13.28]{Washington}), $\Mc_n$ annihilates $A_{n,\psi}$, that means
$\Mc_n\zp[\Ga_n]\subset\Ann_{\zp[\Ga]}A_{n,\psi}$ . Thus
$$
\Mc(T)\La\subset \limp\Ann_{\zp[\Ga]}A_{n,\psi}.
$$

Let $\delta=(\delta_n)_{(n\geq 0)}
\in\limp\Ann_{\zp[\Ga]}A_{n,\psi}$. Then for any $n\geq 0$
$\delta_n A_{n,\psi}=\{0\}$. On the other hand,
$$
e_\psi X=\limp A_{n,\psi}.
$$
Then $\delta e_\psi X=\{0\}$, that implies
$$
\delta\in \Ann_\La e_\psi X=\Mc(T)\La,
$$
that completes the proof. $\square$
\medskip

\medskip
Let $\overline{\W}_n=\W_n\tens\zp$ the $p$-adic adherence of
$\W_n$. The map $\beta$ of Proposition \ref{prop1} induces the map
\begin{eqnarray*}
\overline{\W}_n^-\longrightarrow
(\Ann_{\zz[G_n]}Cl(k_n))^-\tens_\zz\zp=(\Ann_{\zp[G_n]} A_n)^-\\
w\tens a \longmapsto
a\beta(w)\qquad\qquad\qquad\qquad\qquad\qquad\quad
\end{eqnarray*}
that we shall always note $\beta$. Thus we have the short exact
sequence
$$
0\longrightarrow\Hom_{\zz[G_n]}(I_n,\mu_{2p^{n+1}})\tens_{\zz}\zp\longrightarrow
\overline{\W}_n^- \longrightarrow
(\Ann_{\zp[G_n]}A_n)^-\longrightarrow 0.
$$
Applying $e_\psi$ to all the terms of this sequence we get  an
isomorphism of $\zp[G_n]$-modules
\begin{equation}\label{corresp}
e_\psi \overline{\W}_n\simeq e_\psi \Ann_{\zp[G_n]}
A_n=\Ann_{\zp[\Ga_n]}A_{n,\psi}.
\end{equation}
\medskip

Let $z\in \overline{\W}_{n,\psi}=e_\psi \overline{\W}_n$. Then $z$
induces naturally by class field theory (see
\cite[p.455]{Iwasawa}) a morphism of $\zp[\Ga_n]$-modules:
\begin{equation}\label{w_n-map}
z\,:\, \X_{n,\psi}\longrightarrow\U_{n,\psi}.
\end{equation}
\begin{Lem}\label{lem_z}
Let $z\in \overline{\W}_{n,\psi}$ such that $\beta(z)\in
\qq_p[\Ga_n]^*$. Then the kernel of $z$ is
$\Tor_{\zp[\Ga_n]}\X_{n,\psi}$.
\end{Lem}
\noindent{\bf Proof:} We have
$z(\U_{n,\psi})=\U_{n,\psi}^{\beta(z)}$. As $\beta\in
\qq_p[\Ga_n]^*$, the quotient $\U_{n,\psi}/\U_{n,\psi}^{\beta(z)}$
if finite. Thus
$$
\rank_{\zp}z(\X_{n,\psi})=
\rank_{\zp}\X_{n,\psi}=\rank_{\zp}\U_{n,\psi}.
$$
Thus
$\ker(z:\X_{n,\psi}\longrightarrow\U_{n,\psi})=\Tor_{\zp[\Ga_n]}\X_{n,\psi}$.
 $\square$

\medskip

For any $n\in\nn$ let $\Mc_n\in\zp[G_n]$ such that
$$
\Mc_n\equiv \Mc \mod \om_n.
$$
Then $(\Mc_n)_{n\geq0}=\Mc$ in $\La$. Let $w_n\in
\overline{\W}_{n,\psi}=e_\psi\overline{\W}_n$ be the element of
$\overline{\W}_{n,\psi}$ corresponding to $\Mc_n$ via the
homomorphism (\ref{corresp}).
\begin{Rem} The Lemma \ref{lem_z} is applicable to $w_n$, and to the map
that consists in multiplication by $e_\psi\theta_n$.
\end{Rem}
\begin{Lem}\label{ker_w_n}
Let $\overline{J}_n=J_n\tens_{\zz}\zp$ be the $p$-adic adherence
of $J_n$ in $\U_n$ and $\overline{W}_n=W_n\tens_\zz\zp$ the
$p$-adic adherence of $W_n$ in $\U_n$. Then
$$
e_\psi \bar{J}_n\subset  w_n(\X_{n,\psi})\subset e_\psi
\overline{W}_n.
$$
\end{Lem}
\noindent{\bf Proof :} One can verify that
$\beta(e_\psi\overline{\Jc}_n)=e_\psi\theta_n\zp[\Ga_n]$ (see
\cite[Chap. 7]{Washington}). By the Main Conjecture (see \cite[\S
13.6]{Washington})
$$
e_\psi\theta_n\zp[\Ga_n]\subset \Mc_n\zp[\Ga_n].
$$

Set $\widetilde{\W}_{n,\psi}$ the sub-$\zp[\Ga_n]$-module of
$\overline{\W}_{n,\psi}$ generated by $w_n$. Then
$\beta(\widetilde{\W}_{n,\psi})=\Mc_n\zp[\Ga_n]$. Thus
$$
\beta(e_\psi\overline{\Jc}_n)\subset\beta(\widetilde{\W}_{n,\psi}).
$$
As $\beta$ is an isomorphism, this is equivalent to
$$
e_\psi\overline{\Jc}_n\subset\widetilde{\W}_{n,\psi},
$$
That implies
$$
e_\psi\bar{J}_n\subset w_n(\X_{n,\psi}). \ \square
$$

\medskip

Take $(z_n)_{n\geq1}$, $z_n\in \widetilde{W}_{n,\psi}$ such that
$\forall n\geq 1$, $\Res_{n+1,n}\beta(z_{n+1})=\beta(z_n)$. By the
class field theory the following diagram is commutative:
$$
\xymatrix{\X_{n+1,\psi}\ar[r]^{ z_{n+1}}\ar[d]^{\Res_{n+1,n}} & \U_{n+1,\psi}\ar[d]^{N_{n+1,n}}\\
 \X_{n,\psi}\ar[r]^{z_n} & \U_{n,\psi}}
$$
so the map
\begin{equation}
z_\infty\, :\, \X_{\infty,\psi}\longrightarrow\U_{\infty,\psi}
\end{equation}
is naturally well defined and
$$
z_\infty(\X_{\infty,\psi})=\limp z_n(\X_{n,\psi})\subseteq
\U_{\infty,\psi}
$$
\begin{Lem}
The kernel of $z_\infty$ is isomorphic to
$\alpha(e_{\om\psi^{-1}}X)$, where $\alpha(e_{\om\psi^{-1}}X)$ is
the Iwasawa adjoint module of $e_{\om\psi^{-1}}X$.
\end{Lem}
\noindent{\bf Proof :} By the definition of $z_\infty$, $\ker
z_\infty=\limp \ker z_n=\limp \Tor_{\zp[\Ga_n]}\X_{n,\psi}$. But
$\limp \Tor_{\zp[\Ga_n]}\X_{n,\psi}\simeq
\alpha(e_{\omega\psi^{-1}}X)$ (see \cite[Proposition
3.1]{Thong}).$\square$
\medskip

Take $z_n=w_n$ $\forall n\geq 1$. Then $e_\psi\overline
J_\infty\subset w_\infty(\X_{\infty,\psi})$ by the Lemma
\ref{ker_w_n}.

\medskip
\begin{Lem}\label{lem_im_U}
$$
w_\infty(e_\psi\U_\infty)=e_\psi \Mc\U_\infty
$$
\end{Lem}
\noindent{\bf Proof :} Obvious as $e_\psi\U_\infty $ is free of
rank 1. $\square$

\begin{Lem}\label{im_W} The module $W_{\infty,\psi}=\limp \overline{W}_{n,\psi}$
is pseudo-isomorphic to $w_\infty(\X_{\infty,\psi})$.
\end{Lem}
\noindent{\bf Proof:} Let $E$ be the elementary $\La$-module such
that
$$
0\longrightarrow e_\psi X\longrightarrow E\longrightarrow
B\longrightarrow 0,
$$
where $B$ is a finite $\La$-module. Then $\forall n\gg0$,
$\omega_n B=\{0\}$, and by the snake lemma we obtain the exact
sequence
$$
0\longrightarrow B\longrightarrow e_\psi A_n \longrightarrow
E/\omega_n E \longrightarrow B\longrightarrow 0.
$$
Let $Y_n=\Ann_{\zp[\Ga_n]}E/\omega_nE=\Mc_n\zp[\Ga_n]$. It is a
submodule of $Z_n=\Ann_{\zp[\Ga_n]}e_\psi A_n\simeq
e_\psi\overline{\W}_n$, so there exists a submodule
$\widehat{\W}_n$ of $e_\psi\overline{\W}_n$ such that
$\widehat{\W}_n\simeq Y_n$ as $\zp[\Ga_n]$-modules. $Y_n$ being
monogenous, the same is for $\widetilde{\W}_{n,\psi}$, so it is
generated by $w_n$. Thus $\widehat{\W}_n=\widetilde{W}_{n,\psi}$.

There exists $\delta\in \La$, prime to $\Mc$, such that $\delta
B=\{0\}$. Then $\delta Z_n\subset Y_n$, i.e. $\delta
e_\psi\overline{\W}_n\subset \widetilde{\W}_{n,\psi}$. In
particulary that means
\begin{equation}\label{inclusion}
\delta e_\psi\overline{W}_n\subset \widetilde{W}_{n,\psi} \subset
e_\psi \overline{W}_n,
\end{equation}
where $\widetilde{W}_n=w_n(\X_{n,\psi})$.  So, taking the
projective limit in (\ref{inclusion}) we obtain
$$
\delta e_\psi\overline{W}_\infty\subset w_\infty(e_\psi\X_\infty)
\subset e_\psi \overline{W}_\infty.
$$
Thus the quotient module $e_\psi
\overline{W}_\infty/w_\infty(e_\psi\X_\infty)$ is annihilated by
two relatively prime polynomials $\delta$ and $\Mc$, i.e. is
finite (see \cite[\S 13.2]{Washington}). $\square$

\bigskip

The classical class field theory sequence
$$
0\longrightarrow \overline{\Oc_{k_n}^*\cap
\U_n}\longrightarrow\U_n\longrightarrow \X_n\longrightarrow
A_n\longrightarrow 0
$$
gives by taking the $\psi$-parts the short exact sequence
\begin{equation}
0\longrightarrow e_\psi\U_n\longrightarrow e_\psi\X_n
\longrightarrow e_\psi A_n\longrightarrow 0,
\end{equation}
as $\psi\neq \omega$.

Passing to the projective limit in this sequence we obtain the
short exact sequence
\begin{equation}\label{class field sequence}
0\longrightarrow e_\psi\U_\infty\longrightarrow e_\psi\X_\infty
\longrightarrow e_\psi X\longrightarrow 0.
\end{equation}

\begin{Th}
$$
\Char_{\La}(e_\psi\overline{W}_\infty/e_\psi\overline{J}_\infty)=(\Char_\La
e_\psi X_{\infty})/\Mc(T).
$$
\end{Th}
 \noindent{\bf Proof:}
By the Lemma \ref{lem_im_U}, the map $w_\infty$ gives rise to the
map
$$
\overline{w}_\infty\,:\, \
\frac{e_\psi\X_\infty}{e_\psi\U_\infty}\longrightarrow
\frac{e_\psi\U_\infty}{\Mc e_\psi\U_\infty}.
$$
So in virtue of the sequence (\ref{class field sequence}), we have
the map
$$
\overline{w}_\infty\,:\, \ e_\psi X\longrightarrow
\frac{e_\psi\U_\infty}{\Mc e_\psi\U_\infty}.
$$
The kernel of the map $w_\infty$ being the $\La$-torsion module
isomorphic to $\alpha(e_{\omega\psi^{-1}}X)$ and $e_\psi\U_\infty$
being a free $\La$-module, $\ker(w_\infty)\cap
e_\psi\U_\infty=\{0\}$. So
$\ker(\bar{w}_\infty)\simeq\alpha(e_{\omega\psi^{-1}}X)$. And we
have the following exact sequence
\begin{equation*}\label{sequence1}
0\longrightarrow \alpha(e_{\omega\psi^{-1}}X) \longrightarrow
e_\psi X\longrightarrow e_\psi\U_\infty/\Mc e_\psi\U_\infty.
\end{equation*}
Let $F=\Char_\La e_\psi X$ the characteristic polynomial of
$e_\psi X$ and $\theta_{\infty,\psi}=(e_\psi\theta_n)_{n\geq 0}$.
Then by the Main Conjecture  and by the Lemma \ref{ker_w_n} we
obtain the second exact sequence
\begin{equation*}\label{sequence2}
0\longrightarrow \alpha(e_{\omega\psi^{-1}}X) \longrightarrow
e_\psi X\stackrel{\times\theta_{\infty,\psi}}\longrightarrow
e_\psi\U_\infty/F e_\psi\U_\infty
\end{equation*}
 These two sequences give two short exact
sequences
\begin{equation*}\label{sequence3}
0\longrightarrow \alpha(e_{\omega\psi^{-1}}X) \longrightarrow
e_\psi X\longrightarrow w_\infty(\X_{\infty,\psi})/\Mc
e_\psi\U_\infty\longrightarrow 0
\end{equation*}
and
\begin{equation*}\label{sequence4}
0\longrightarrow \alpha(e_{\omega\psi^{-1}}X) \longrightarrow
e_\psi X\longrightarrow e_\psi\overline{J}_\infty/F
e_\psi\U_\infty\longrightarrow 0,
\end{equation*}
as $\theta_\infty\X_{\infty,\psi}=e_\psi\overline{J}_\infty$ (see
\cite[Lemme 8]{moi}). Thus
\begin{equation}\label{first_equation}
\Char_\La  w_\infty(\X_{\infty,\psi})/\Mc
e_\psi\U_\infty=\Char_\La e_\psi\overline{J}_\infty/F
e_\psi\U_\infty.
\end{equation}
Set $e_\psi\widetilde{W}_\infty=w_\infty(\X_{\infty,\psi})$. The
tautological short exact sequence
$$
0\longrightarrow e_\psi \overline{J}_{\infty} \longrightarrow
e_\psi\widetilde{W}_\infty\longrightarrow
e_\psi\widetilde{W}_\infty/e_\psi \overline{J}_{\infty}
\longrightarrow 0
$$
gives rise to the short exact sequence
\begin{equation*}\label{first_sequence}
0\longrightarrow \frac{e_\psi
\overline{J}_{\infty}}{F\U_{\infty\psi}} \longrightarrow
\frac{e_\psi\widetilde{W}_\infty}{F\U_{\infty,\psi}}\longrightarrow
\frac{e_\psi\widetilde{W}_\infty}{e_\psi \overline{J}_{\infty}}
\longrightarrow 0.
\end{equation*}
Thus
\begin{equation}\label{second_equation}
    \Char\frac{e_\psi\overline{J}_{\infty}}{F\U_{\infty,\psi}}=
    \Char\frac{e_\psi\widetilde{W}_\infty}{F\U_{\infty,\psi}}
    \left(\Char\frac{e_\psi\widetilde{W}_\infty}{e_\psi
    \overline{J}_{\infty}}\right)^{-1}.
\end{equation}
In the same way, the sequence
$$
0\longrightarrow \Mc e_\psi\U_\infty \longrightarrow
e_\psi\widetilde{W}_\infty\longrightarrow
e_\psi\widetilde{W}_\infty/\Mc e_\psi\U_\infty \longrightarrow 0
$$
gives rise to the sequence
$$
0\longrightarrow \frac{\Mc \U_{\infty,\psi}
}{F\U_{\infty\psi}}\longrightarrow
\frac{e_\psi\widetilde{W}_\infty}{F\U_{\infty\psi}}\longrightarrow
\frac{e_\psi\widetilde{W}_\infty}{\Mc\U_{\infty,\psi}}
\longrightarrow 0
$$
Thus
\begin{equation}\label{third_equation}
    \Char\frac{e_\psi\widetilde{W}_\infty}{\Mc\U_{\infty,\psi}}=
    \Char\frac{e_\psi\widetilde{W}_\infty}{F\U_{\infty\psi}}\left(\Char\frac{\Mc \U_{\infty,\psi}
}{F\U_{\infty\psi}}\right)^{-1}.
\end{equation}
 Comparing the equalities
(\ref{first_equation}), (\ref{second_equation}) and
(\ref{third_equation}) we obtain the equality
$$
\Char\frac{\Mc \U_{\infty,\psi}}{F\U_{\infty,\psi}}=
\Char\frac{e_\psi\widetilde{W}_\infty}{e_\psi
    \overline{J}_{\infty}}.
$$
In virtue of the Lemma \ref{im_W},
$$
\Char\frac{e_\psi\widetilde{W}_\infty}{e_\psi
    \overline{J}_{\infty}}=\Char\frac{e_\psi\overline{W}_\infty}{e_\psi
    \overline{J}_{\infty}}.
$$
As $\U_{\infty,\psi}$ is free of rank 1,
$$
\Char\frac{\Mc
\U_{\infty,\psi}}{F\U_{\infty,\psi}}=\frac{F(T)}{\Mc(T)}.
$$
So
$$
\Char\frac{e_\psi\overline{W}_\infty}{e_\psi
    \overline{J}_{\infty}}=\frac{F(T)}{\Mc(T)}. \ \square
$$
\begin{Cor}[cf. {\cite[Th\'eor\`eme 1]{moi}}]~\label{Cor1}
$$
\Char_\La \frac {\U_{\infty,\psi}}{e_\psi \overline{J}\infty} =
\Char_\La \alpha(e_{\omega\psi^{-1}}X).
$$
\end{Cor}

\begin{Cor} The module $e_\psi X$ is pseudo-monogenous if and only
if  the quotient module
$e_\psi\overline{W}\infty/e_\psi\overline{J}\infty$ is finite.
\end{Cor}

By the corollary \ref{Cor1}, we see that  Greenberg Conjecture
implies that
$$
\frac{e_\psi\U_{\infty}}{e_\psi \overline{W}_\infty} \text{ is
finite.}
$$
So it is natural to ask the following question.

\noindent{\bf Question:} Let $p$ be an odd prime number. Let
$\psi$ be an odd character of $\Gal(\qq(\z_p)/\qq)$,
$\psi\neq\omega$, then, is it true that
$$
\frac{e_\psi\U_{\infty}}{e_\psi \overline{W}_\infty} \text{ is
finite ?}
$$
\begin{Rem}
Note that the positive answer to this question is equivalent to
$$
\Char e_\psi X=\Char \alpha(e_{\omega\psi^{-1}}X)\times \Mc,
$$
so it implies weak Greenberg Conjecture (see \cite{Raphi_Thong},
\cite{Thong2}).
\end{Rem}

\bigskip

\noindent
Tatiana BELIAEVA, Bruno ANGLES\\
Universit\'e de Caen,\\
 Laboratoire Nicolas Oresme, CNRS 6139,\\
 Campus II, Boulevard Mar\'echal Juin,\\
 B.P. 5186, 14032 Caen Cedex,, France.\\

\noindent
 bruno.angles@math.unicaen.fr\\ tatiana.beliaeva@math.unicaen.fr\\

\end{document}

%% file: def_engl.tex
\newtheorem{Th}{Theorem}

\newtheorem{Lem}{Lemma}
\newtheorem{Prop}{Proposition}
\newtheorem{Cor}{Corollary}
\newtheorem{Def}{Definition}
\newtheorem{Rem}{Remark}
\newcommand{\zz}{\mathbb Z}
\newcommand{\nn}{\mathbb N}

\newcommand{\qq}{\mathbb Q}

\newcommand{\limp}{\displaystyle \lim_{\longleftarrow}}
\newcommand{\g}{\gamma}

\newcommand{\Ga}{\Gamma}
\newcommand{\La}{\Lambda}
\newcommand{\D}{\Delta}
\newcommand{\om}{\omega}
\newcommand{\U}{{\cal U}}

\newcommand{\X}{{\mathfrak X}}

\newcommand{\ag}{{\mathfrak a}}

\newcommand{\lc}{{\mathfrak l}}
\newcommand{\eps}{\varepsilon}
\newcommand{\z}{\zeta}

\newcommand{\zp}{\zz_p}
\newcommand{\qp}{\qq_p}

\newcommand{\tens}{\otimes}

\newcommand{\Oc}{{\cal O}}
\newcommand{\W}{{\cal W}}
\newcommand{\Mc}{{\cal M}}
\newcommand{\Jc}{{\cal J}}

\newcommand{\Hom}{{\rm Hom\,}}
\newcommand{\Res}{{\rm Res\,}}
\newcommand{\Ann}{{\rm Ann\,}}

\newcommand{\Gal}{{\rm Gal\,}}

\newcommand{\Tr}{{\rm Tr}\,}
\newcommand{\Char}{{\rm char}\,}

\newcommand{\Tor}{{\rm Tor}\,}

\newcommand{\rank}{{\rm rank}\,}
\newcommand{\Lc}{{\mathfrak L}}
\newcommand{\p}{{\mathfrak p}}
